\documentclass[12pt]{article}

\def\sqr#1#2{{\vcenter{\vbox{\hrule height.#2pt
              \hbox{\vrule width.#2pt height#1pt \kern#1pt \vrule width.#2pt}
              \hrule height.#2pt}}}}
\def\signed #1{{\unskip\nobreak\hfil\penalty50
              \hskip2em\hbox{}\nobreak\hfil#1
              \parfillskip=0pt \finalhyphendemerits=0 \par}}
\def\endpf{\signed {$\sqr69$}}

\def\dbR{{\mathop{\rm l\negthinspace R}}}
\def\dbC{{\mathop{\rm l\negthinspace\negthinspace\negthinspace C}}}

\def\3n{\negthinspace \negthinspace \negthinspace }
\def\2n{\negthinspace \negthinspace }
\def\1n{\negthinspace }

\def\dbC{{\mathop{\rm l\negthinspace\negthinspace\negthinspace C}}}

\def\ds{\displaystyle}

\def\dbN{{\mathop{\rm l\negthinspace N}}}

\def\dbR{{\mathop{\rm l\negthinspace R}}}
\def\={\buildrel \triangle \over =}

\def\a{\alpha}
\def\b{\beta}
\def\g{\gamma}
\def\d{\delta}

\def\l{\lambda}

 \def\n{\nabla}

\def\t{\times}

\def\th{\theta}
\def\o{\omega}

\def\i{\infty}
\def\ns{\noalign{\ss} }

\def\G{\Gamma}
\def\D{\Delta}

\def\Si{\Sigma}

\def\O{\Omega}

\def\cC{{\cal C}}

\def\cP{{\cal P}}

\def\cl{{\cal l}}
\def\no{\noindent}

\def\ms{\medskip}
\def\bs{\bigskip}
\def\q{\quad}
\def\qq{\qquad}
\def\hb{\hbox}

\def\min{\mathop{\rm min}}

\def\pa{\partial}

\def\cd{\cdot}

\def\div{\hbox{\rm div$\,$}}

\def\cl{\overline}

\def\|{\Big |}
\def\({\Big (}
\def\){\Big )}
\def\[{\Big[}
\def\]{\Big]}
\def\be{\begin{equation}}
\def\bel{\begin{equation}\label}
\def\ee{\end{equation}}
\def\bt{\begin{theorem}}
\def\bcd{\begin{condition}}
\def\ecd{\end{condition}}
\def\et{\end{theorem}}
\def\bc{\begin{corollary}}
\def\ec{\end{corollary}}
\def\bde{\begin{definition}}
\def\ede{\end{definition}}
\def\bl{\begin{lemma}}
\def\el{\end{lemma}}
\def\bp{\begin{proposition}}
\def\ep{\end{proposition}}
\def\br{\begin{remark}}
\def\er{\end{remark}}
\def\ba{\begin{array}}
\def\ea{\end{array}}
\def\ed{\end{document}}
\def\ns{\noalign{\ms}}
\def\ds{\displaystyle}

\def\square#1{\vbox{\hrule\hbox{\vrule height#1%
     \kern#1\vrule}\hrule}}
\def\rectangle#1#2{\vbox{\hrule\hbox{\vrule height#1%
     \kern#2\vrule}\hrule}}

\font\tenbb=msbm10 \font\sevenbb=msbm7 \font\fivebb=msbm5

\newfam\bbfam
\scriptscriptfont\bbfam=\fivebb \textfont\bbfam=\tenbb
\scriptfont\bbfam=\sevenbb

\def\oO{{\overline \O}}
\def\ov{{\overline v}}

\newtheorem{lemma}{Lemma}[section]
\newtheorem{remark}{Remark}[section]

\newtheorem{theorem}{Theorem}[section]
\newtheorem{corollary}{Corollary}[section]

\newtheorem{definition}{Definition}[section]
\newtheorem{proposition}{Proposition}[section]
\newtheorem{condition}{Condition}[section]

\makeatletter
   
   \@addtoreset{equation}{section}
\makeatother

\begin{document}
\title{\bf Null controllability for the parabolic equation with a complex principal
part\thanks{This work is partially supported by the NSF of China
under grants 10525105 and 10771149. The author gratefully
acknowledges Professor  Xu Zhang for his guidance, encouragement and
suggestions.}}

\author{Xiaoyu Fu\thanks{School of Mathematics,  Sichuan University,
Chengdu 610064,  China. {\small\it E-mail:} {\small\tt
rj\_xy@163.com}.}}

\date{}

\maketitle

 \begin{abstract}
\no This paper is addressed to a study of the null controllability
for the semilinear parabolic equation with a complex principal part.
For this purpose, we establish a key weighted identity for partial
differential operators
$(\a+i\b)\pa_t+\sum\limits_{j,k=1}^n\pa_k(a^{jk}\pa_j)$ (with real
functions $\a$ and $\b$), by which we develop a universal approach,
based on global Carleman estimate, to deduce not only the desired
explicit observability estimate for the linearized complex
Ginzburg-Landau equation, but also all the known
controllability/observability results for the parabolic, hyperbolic,
Schr\"odinger and plate equations that are derived via Carleman
estimates.
 \end{abstract}

\bs

\no{\bf 2000 Mathematics Subject Classification}.  Primary 93B05;
Secondary 35B37, 35L70, 93B07, 93C20.

\bs

\no{\bf Key Words}. Null controllability, observability, universal
approach, Carleman estimate,   Parabolic equation with a complex
principal part.

\section{Introduction and main results}\label{is1}

Given $T>0$ and a bounded domain $\O$ of $\dbR^{n}$ ($n\in\dbN$)
with $C^2$ boundary $\G$. Fix an open non-empty subset $\o$ of $\O$
and denote by $\chi_\o$ the characteristic function of $\o$. Let
$\o_0$ be another non-empty open subset of $\O$ such that $\overline
{\o_0}\subset\o$.  Put
$$Q=(0,T)\times\O,\q\Si=(0,T)\times\G,\q Q_0=(0,T)\times\o_0.$$
In the sequel, we will use the notation $ y_j= y_{ x_j}$, where
$x_j$ is the $j$-th coordinate of a generic point
$x=(x_1,\cdots,x_n)$ in $\dbR^n$. In a similar manner, we use the
notation $ z_j$, $v_j$, etc. for the partial derivatives of $z$ and
$v$ with respect to $x_j$. Throughout this paper, we will use
$C=C(T,\O,\o)$ to denote generic positive constants which may vary
from line to line (unless otherwise stated). For any $c\in\dbC$, we
denote its complex conjugate by $\overline c$.

Fix $a^{jk}(\cd)\in C^{1,2}(\cl{Q}; \;\dbR)$ satisfying
 \bel{a1}
a^{jk}(t,x)=a^{kj}(t,x),\q (t,x)\in\cl{Q},~j,k=1,2,\cdots,n,
 \ee
and for some constant $s_0>0$,
 \bel{a2}
  \sum_{j,k}a^{jk}\xi_j\overline{\xi_k}\ge s_0|\xi|^2,\qq (t,x,\xi)\equiv
 (t,x,\xi_1,\cdots,\xi_n)\in\overline Q\t\dbC^n.
 \ee
Next, we fix a function $f(\cd)\in C^1(\,\dbC)$ satisfying $f(0)=0$
and
 the following condition:
 \bel{fn1}
  \lim_{s\to\i} \frac{|f(s)|}{|s|\ln^{1/2} |s|}=0.
 \ee
Note that $f(\cd)$ in the above can have a superlinear growth. We
are interested in the following semilinear parabolic equation with a
complex principal part:
 \bel{3ga1}\left\{\ba{ll}\ds (1+ib)y_t-\sum_{j,k=1}^n(a^{jk}y_j)_k=\chi_\o u+f(y) \q &\mbox {  in $Q$},\\
 \ns\ds y=0 &\mbox { on $\Si$},\\
 \ns\ds  y(0)=y_0& \mbox { in $\O$},
 \ea\right.\ee
where $i=\sqrt{-1}$, $b\in \dbR$. In (\ref{3ga1}), $y=y(t,x)$ is the
state and $u=u(t,x)$ is the control.

One of our main objects in this paper is to study the null
controllability of system (\ref{3ga1}), by which we mean that, for
any given initial state $y_0$, find (if possible) a control $u$ such
that the weak solution of (\ref{3ga1}) satisfies $y(T)=0$.

In order to derive the null controllability of (\ref{3ga1}), by
means of the well-known duality argument (see \cite[p.282, Lemma
2.4]{LY}, for example), one needs to consider the following dual
system of the linearized system of (\ref{3ga1}) (which can be
regarded as a linearized complex Ginzburg-Landau equation):
   \bel{gap1}\left\{\ba{ll}\ds
Gz=q(t,x)z &\mbox { in }
 Q,\\
 \ns\ds z=0 &\mbox { on }
 \Si,\\
\ns\ds z(T)=z_T &\mbox { in } \O,
 \ea\right.\ee
 where $q(\cd)\in L^\i(0,T;L^n(\O))$ is a potential and
 \bel{cp2}
Gz\=(1+ib)z_t+\sum_{j,k=1}^n(a^{jk}z_j)_k.
 \ee

By means of global Carleman inequality, we shall establish the
following explicit observability estimate for solutions of system
(\ref{gap1}).

 \bt\label{gt5}  Let $a^{jk}(\cd)\in C^{1,2}(\cl{Q}; \;\dbR)$ satisfy (\ref{a1})--(\ref{a2}), $q(\cd)\in L^\i(0,T;L^n(\O))$ and $b\in\dbR$.  Then
 there is a constant $C>0$ such that for all solutions of system (\ref{gap1}), it holds
  \bel{3g1}
 |z(0)|_{L^2(\O)}\le \cC(r)|z|_{L^2((0,T)\t\o)},\qq\forall\; z_T\in
 L^2(\O)
  \ee
 where
 \bel{ccg}
 \cC(r)\= C_0e^{C_0r^2}, \q C_0=C(1+b^2
),\q r\= |q|_{L^\i(0,T;L^n(\O))}.
 \ee
 \et

Thanks to the dual argument and the fixed point technique, Theorem
\ref{gt5} implies the following controllability result for system
(\ref{3ga1}).

\bt\label{gt11}  Let $a^{jk}(\cd)\in C^{1,2}(\cl{Q}; \;\dbR)$
satisfy (\ref{a1})--(\ref{a2}), $f(\cd)\in C^1(\,\dbC)$ satisfy
$f(0)=0$ and (\ref{fn1}), and $b\in\dbR$. Then for any given $y_0\in
L^2(\O)$, there is a control $u\in L^2((0,T)\t\o)$ such that the
weak solution $y(\cdot)\in C([0,T];L^2(\O))\bigcap
L^2(0,T;H_0^1(\O))$ of system (\ref{3ga1}) satisfies $y(T)=0$ in
$\O$.
 \et

The controllability problem for system (\ref{3ga1}) with $b=0$
(i.e., linear and semilinear parabolic equations) has been studied
by many authors and it is now well-understood. Among them, let us
mention \cite{DFGZ, DZZ, FI, FZ2} on what concerns null
controllability, \cite{FPZ, FZ, FZ2, Z2} for approximate
controllability, and especially \cite{Zu} for recent survey in this
respect. However, for the case $b\not=0$,  very little is know in
the previous literature. To the best of our knowledge, \cite{Fu} is
the only paper addressing the global controllability for
multidimensional system (\ref{3ga1}). We refer to \cite{RZ1} for a
recent interesting result on local controllability of semilinear
complex Ginzburg-Landau equation.

We remark that condition (\ref{fn1}) is not sharp. Indeed, following
\cite{DFGZ}, one can establish the null controllability of system
(\ref{3ga1}) when the nonlinearity $f(y)$ is replaced by a more
general form of $f(y,\nabla y)$ under the assumptions that
$f(\cd,\cd)\in C^1(\,\dbC^{1+n})$, $f(0,0)=0$ and
 \bel{fssn1}
 \ba{ll}\ds
 \lim_{|(s,p)|\to\i} \frac{\ds\left|\int_0^1 f_s(\tau s,\tau p)d\tau \right|}{\ln^{3/2}(1+ |s|+
 |p|)}=0,\\
  \ns
  \ds
 \lim_{|(s,p)|\to\i} \frac{\ds\left|\left(\int_0^1 f_{p_1}(\tau s,\tau p)d\tau,\cdots\int_0^1 f_{p_n}(\tau s,\tau p)d\tau \right)
 \right|}{\ln^{1/2}(1+ |s|+ |p|)}=0,
 \ea
 \ee
where $p=(p_1,\cdots,p_n)$. Moreover, following \cite{DZZ}, one can
show that the assumptions on the growth of the nonlinearity
$f(y,\nabla y)$ in (\ref{fssn1}) are sharp in some sense. Since the
techniques are very similar to \cite{DFGZ, DZZ}, we shall not give
the details here.

Instead, as a byproduct of the fundamental identity  established in
this work (to show the observability inequality (\ref{3g1})), we
shall develop a universal approach for controllability/observability
problems governed by partial differential equations (PDEs for
short), which is the second main object of this paper. The study of
controllability/observability problem for PDEs  began in the 1960s,
for which various techniques have been developed in the last decades
(\cite{BLR, Coron, FI, Lions, Zu}). It is well-known that the
controllability of PDEs depends strongly on the nature of the
system, say time reversibility or not. Typical examples are the wave
and heat equations. It is clear now that there exists essential
differences between the controllability/observability theories for
these two equations. Naturally, one expects to know whether there
are some relationship between these two systems of different nature.
Especially, it would be quite interesting to establish a unified
controllability/observability theory for parabolic and hyperbolic
equations. This problem was posed by D.L. Russell in \cite{Ru},
where one can also find some preliminary result; further results are
available in \cite{LZZ, Zh}. In \cite{LW}, the authors analyzed the
controllability/observability problems for PDEs from the point of
view of methodology. It is well known that these problems may be
reduced to the obtention of suitable observability inequalities for
the underlying homogeneous systems. However, the techniques that
have been developed to obtain such estimates depend heavily on the
nature of the equations. In the context of the wave equation one may
use multipliers (\cite{Lions}) or microlocal analysis (\cite{BLR});
while, in the context of heat equations, one uses Carleman estimates
(\cite{FZ, FI}). Carleman estimates can also be used to obtain
observability inequalities for the wave equation (\cite{Zh0}).
However, the Carleman estimate that has been developed up to now to
establish observability inequalities of PDEs depend heavily on the
nature of the equations, and therefore a unified Carleman estimate
for those two equations has not been developed before. In this
paper, we present a point-wise weighted identity for partial
differential operators
$(\a+i\b)\pa_t+\sum\limits_{j,k=1}^n\pa_k(a^{jk}\pa_j)$ (with real
functions $\a$ and $\b$), by which we develop a unified approach,
based on global Carleman estimate, to deduce not only the
controllability/observability results for systems (\ref{3ga1}) and
(\ref{gap1}), but also all the known controllability/observability
results for the parabolic, hyperbolic, Schr\"odinger and plate
equations that are derived via Carleman estimates (see Section
\ref{ss2} for more details). We point out that this identity has
other interesting applications, say, in \cite{PZ} it is applied to
derive an asymptotic formula of reconstructing the initial state for
a Kirchhoff plate equation with a logarithmical convergence rate for
smooth data; while in \cite{Fu1} it is applied to establish sharp
logarithmic decay rate for general hyperbolic equations with damping
localized in arbitrarily small set by means of an approach which is
different from that in \cite{B}.

The rest of this paper is organized as follows. In Section
\ref{ss2}, we establish a fundamental point-wise weighted identity
for partial differential operators of second order and give some of
its applications. In Section \ref{ss3}, we derive a modified
point-wise inequality for the parabolic operator. This estimate will
play a key role when we establish in Section \ref{ss4} a global
Carleman estimate for the parabolic equation with a complex
principal part. Finally, we will prove our main results in Section
\ref{ss5}.

\section{A weighted identity for partial differential operators and its applications}\label{ss2}

In this section, we will establish a point-wise weighted identity
for partial differential operators of second order with a complex
principal part, which has an independent interest. First, we
introduce the following second order operator:
  \bel{cp}\ba{ll}\ds
 \cP z\=(\a+i\b)z_t+\sum_{j,k=1}^m(a^{jk}z_j)_k,\qq m\in\dbN.
 \ea\ee
We have the following fundamental identity.

 \bt\label{2t1} Let $\a,\;\b\in C^{2}(\dbR^{1+m};\dbR)$ and $a^{jk}\in
C^{1,2}(\dbR^{1+m};\dbR)$ satisfy $a^{jk}=a^{kj}$
$(j,k=1,2,\cdots,m)$. Let $z,\;v\in C^2(\dbR^{1+m}; \;\dbC)$ and
$\Psi,\;\ell\in C^2(\dbR^{1+m};\dbR)$. Set $\th=e^\ell$ and $v=\th
z$. Then
 \bel{2a2}\ba{ll}\ds
 \th(\cP z\overline {I_1}+\overline{\cP z} I_1)+M_t+\div V\\
\ns\ds =2|I_1|^2+
\sum_{j,k,j',k'=1}^m\[2(a^{j'k}\ell_{j'})_{k'}a^{jk'}-(a^{jk}a^{j'k'}\ell_{j'})_{k'}+\frac{1}{2}(\a
a^{jk})_t-a^{jk}\Psi\](v_k\ov_j+\ov_k
v_j)\\
\ns\ds\q+i\sum_{j,k=1}^m\[(\b a^{jk}\ell_j)_t+a^{jk}(\b\ell_t)_j\](\ov_kv-v_k\ov)-\sum_{j,k=1}^ma^{jk}\a_k(v_j\ov_t+\ov_jv_t)\\
\ns\ds\q+i\[\b\Psi+\sum_{j,k=1}^m(\b a^{jk}\ell_j)_k\](\ov
v_t-v\ov_t)+B|v|^2,
 \ea\ee
 where
 \bel{1f3}\left\{\ba{ll}\ds
 A\=\sum_{j,k=1}^ma^{jk}\ell_j\ell_k-\sum_{j,k=1}^m(a^{jk}\ell_j)_k -\Psi,\\
 \ns\ds
I_1\=i\b v_t-\a\ell_tv+\sum_{j,k=1}^m(a^{jk}v_j)_k+Av,
 \ea\right.\ee
 and
 \bel{1a3}\left\{\ba{ll}\ds
M\=\[(\a^2+\b^2)\ell_t-\a A\] |v|^2+\a\sum_{j,k=1}^ma^{jk}v_j\ov_k+i\b\sum_{j,k=1}^ma^{jk}\ell_j(\ov_kv-v_k\ov),\\
\ns\ds V\=[V^1,\cdots,V^k,\cdots,V^m],\\
 \ns\ds V^k\=\sum_{j,j',k'=1}^m\Big\{-i \b\[a^{jk}\ell_j(\ov_tv-\ov
 v_t)+a^{jk}\ell_t(v_j\ov-\ov_jv)\]-\a a^{jk}(v_j\ov_t+\ov_jv_t)\\
\ns\ds\qq\q +\(2a^{jk'}a^{j'k}-a^{jk}a^{j'k'}\)\ell_j(v_{j'}\ov_{k'}+\ov_{j'}v_{k'}) -\Psi a^{jk}(v_j\ov+\ov_jv)\\
 \ns\ds\qq\q
+a^{jk}(2A\ell_j+\Psi_j-2\a\ell_j\ell_t)|v|^2\Big\},\\
\ns\ds
B\=(\a^2\ell_t)_t+(\b^2\ell_{t})_t-(\a A)_t-2\[\sum_{j,k=1}^m(\a a^{jk}\ell_j\ell_{t})_k+\a\Psi\ell_t\]\\
\ns\ds\qq\q+\sum_{j,k=1}^m(a^{jk}\Psi_k)_j+2\[\sum_{j,k=1}^m(a^{jk}\ell_jA)_k+A\Psi\].
 \ea\right.\ee \et

Several remarks are in order.

 \br We see that only the symmetry condition of $(a^{jk})_{m\t m}$ is
assumed in the above. Therefore, Theorem \ref{2t1} is applicable to
ultra-hyperbolic or ultra-parabolic differential operators.
 \er
 \br
Note that when $\ds\Psi=-\sum_{j,k=1}^m(a^{jk}\ell_j)_k$ and
$\a(t,x)\equiv a,\;\b(t,x)\equiv b$ with $a,\;b\in\dbR$, Theorem
\ref{2t1} is reduced to \cite[Theorem 1.1]{Fu}. Here, we add an
auxiliary function $\Psi$ in the right-hand side of the multiplier
$I_1$ so that the corresponding identity coincides with
\cite[Theorem 4.1]{FYZ} for the case of hyperbolic operators.
Moreover, we will see that the modified identity is better than
\cite{Fu} in some sense.
 \er
 \br\label{rekh}
By choosing $\a(t,x)\equiv1,\; \b(t,x)\equiv0$ and
$\ds\Psi=-2\sum_{j,k=1}^m(a^{jk}\ell_j)_k$ in Theorem \ref{2t1}, one
obtains a weighted identity for the parabolic operator. By this and
following \cite{TZ}, one may recover all the
controllability/observability results for the parabolic equations in
\cite{DFGZ} and \cite{FI}.
 \er
 \br
By choosing  $a^{jk}(t,x)\equiv a^{jk}(x)$ and
$\a(t,x)=\b(t,x)\equiv0$ in Theorem \ref{2t1}, one obtains the key
identity derived in \cite{FYZ} for the controllability/observability
results on the general hyperbolic equations.
 \er
 \br
By choosing  $(a^{jk})_{1\le j,k\le m}$ to be the identity matrix,
$\a(t,x)\equiv0,\;\b(t,x)\equiv1$ and $\Psi=-\D\ell$ in Theorem
\ref{2t1}, one obtains the pointwise identity derived in \cite{LTZ}
for the observability results for the nonconservative Schr\"odinger
equations. Also, this yields the controllability/observability
results in \cite{Z} for the plate equations and the results for
inverse problem for the Schr\"odinger equation in \cite{BP}.
 \er
 \br
By choosing $(a^{jk})_{1\le j,k\le m}$ to be the identity matrix,
$\a(t,x)\equiv 0,\;\b(t,x)\equiv p(x)$ and $\Psi=-\D\ell$ in Theorem
\ref{2t1}, one obtains the pointwise identity for the Schr\"odinger
operator: $\ds ip(x)\pa_t+\D$. Further, by choosing
 $$
\ell(t,x)=s\varphi,\q\varphi=e^{\g(|x-x_0|^2-c|t-t_0|^2)}
 $$
with $\g>0$, $c>0$, $x_0\in \dbR^n\setminus\oO$ and assuming $\n\log
p\cdot(x-x_0)>-2$. One can recover the fundamental Carleman estimate
for Schr\"odinger operator $\ds ip(x)\pa_t+\D$ derived in
\cite[Lemma 2.1]{YY}.
 \er

 {\it Proof of Theorem
\ref{2t1}.} The proof is divided into several steps.

{\it Step 1.} By $\th=e^\ell,\; v=\th z$,  we have (recalling
(\ref{1a3}) for the definition of $I_1$)
 \bel{1a6}\ba{ll}\ds
 \th\cP z&\ds=(\a+i\b)v_t-(\a+i\b)\ell_tv+\sum_{j,k=1}^m(a^{jk}v_j)_k\\
 \ns&\ds\qq\q+\sum_{j,k=1}^ma^{jk}\ell_j\ell_kv-2\sum_{j,k=1}^ma^{jk}\ell_jv_k-\sum_{j,k=1}^m(a^{jk}\ell_j)_kv\\
 \ns&\ds=I_1+I_2,
 \ea\ee
where
 \bel{1a7}
  I_2\=\a v_t-i\b\ell_tv-2\sum_{j,k=1}^ma^{jk}\ell_jv_k+\Psi v.
 \ee
Hence, by recalling (\ref{1f3}) for the definition of $I_1$, we have
 \bel{1a8}
 \th(\cP z\overline {I_1}+\overline{\cP z} I_1)= 2|I_1|^2+(I_1\overline I_2+I_2\overline I_1).
 \ee

\ms

{\it Step 2.} Let us compute $I_1\overline I_2+I_2\overline I_1$.
 Denote the four terms in the right hand side of $I_1$ and $I_2$ by $I_1^j$ and $I_2^j$, respectively, $j=1, 2, 3,
 4$. Then
 \bel{I}
 I_1^1(\overline {I_2^1}+ \overline {I_2^2})+\overline {I_1^1}
 (I_2^1+I_2^2)=-(\b^2\ell_t|v|^2)_t+(\b^2\ell_{t})_t|v|^2.
 \ee
Note that
 $$\left\{\ba{ll}\ds
2v\ov_t=(|v|^2)_t-(\ov v_t-v\ov_t),\\
\ns \ds2v\ov_k=(|v|^2)_k-(\ov v_k-v\ov_k).
 \ea\right.$$
Hence, we get
 \bel{i13}\ba{ll}
 I_1^1 (\overline {I_2^3}+\overline {I_2^4})+\overline
 {I_1^1}(I_2^3+I_2^4)\ds\\
 \ns\ds=-2i\sum_{j,k=1}^m\[(\b a^{jk}\ell_jv\ov_k)_t-(\b a^{jk}\ell_j)_tv\ov_k\]\\
 \ns\ds\q+2i\sum_{j,k=1}^m\[(\b a^{jk}\ell_jv\ov_t)_k-(\b a^{jk}\ell_j)_kv\ov_t\]-i\b\Psi(v\ov_t-v_t\ov)\\
\ns\ds=-i\sum_{j,k=1}^m\[\b a^{jk}\ell_j(v\ov_k-\ov v_k)\]_t+i\sum_{j,k=1}^m\[\b a^{jk}\ell_j(v\ov_t-\ov v_t)\]_k\\
\ns\ds\q-i\sum_{j,k=1}^m(\b a^{jk}\ell_j)_t(\ov
v_k-v\ov_k)+i\[\b\Psi+\sum_{j,k=1}^m(\b a^{jk}\ell_j)_k\](\ov
v_t-v\ov_t).
 \ea\ee
Next,
 \bel{i12i2}\ba{ll}&\ds
I_1^2(\overline I_2^1+\overline I_2^2+\overline I_2^3+\overline
I_2^4)+\overline I_1^2
(I_2^1+I_2^2+I_2^3+I_2^4)\\
 \ns&\ds=-(\a^2\ell_t|v|^2)_t+2\sum_{j,k=1}^m(\a a^{jk}\ell_j\ell_t|v|^2)_k\\
\ns&\ds\q+(\a^2\ell_{t})_t|v|^2-2\[\sum_{j,k=1}^m(\a
a^{jk}\ell_j\ell_{t})_k+\a\Psi\ell_t\]|v|^2.
 \ea\ee
Noting that $a^{jk}=a^{kj}$, we have
 \bel{i31}\ba{ll}\ds
I_1^3 (\overline {I_2^1}+\overline {I_2^2})+\overline
 {I_1^3}(I_2^1+I_2^2)\\
\ns\ds
=\sum_{j,k=1}^m\[\a a^{jk}(v_j\ov_t+\ov_jv_t)\]_k-\sum_{j,k=1}^ma^{jk}\a_k(v_j\ov_t+\ov_jv_t)-\sum_{j,k=1}^m(\a a^{jk}v_j\ov_k)_t\\
\ns\ds\q +\sum_{j,k=1}^m(\a a^{jk})_tv_j\ov_k+i\sum_{j,k=1}^m\[\b
a^{jk}\ell_t(v_j\ov-\ov_j
 v)\]_k+i\sum_{j,k=1}^ma^{jk}(\b\ell_{t})_k(\ov_j
 v-v_j\ov)\\
 \ns\ds=\sum_{j,k=1}^m\[\a a^{jk}(v_j\ov_t+\ov_jv_t)+i\b a^{jk}\ell_t(v_j\ov-\ov_j v)\]_k\\
 \ns\ds\q-\sum_{j,k=1}^m(\a a^{jk}v_j\ov_k)_t-\sum_{j,k=1}^ma^{jk}\a_k(v_j\ov_t+\ov_jv_t)\\
 \ns\ds\q+{1\over2}\sum_{j,k=1}^m(\a a^{jk})_t(v_j\ov_k+v_k\ov_j)+i\sum_{j,k=1}^ma^{jk}(\b\ell_{t})_k(\ov_j
 v-v_j\ov).
 \ea\ee
Using the symmetry condition of $a^{jk}$ again, we obtain
  \bel{ii33}\ba{ll}&\ds
2\sum_{j,k,j',k'=1}^ma^{jk}a^{j'k'}\ell_j(v_{j'}\ov_{kk'}+\ov_{j'}v_{kk'})\\
 \ns&\ds=\sum_{j,k,j',k'=1}^m\[a^{jk}a^{j'k'}\ell_j(v_{j'}\ov_{k'}+\ov_{j'}v_{k'})\]_k-\sum_{j,k,j',k'=1}^m(a^{jk}a^{j'k'}\ell_j)_k(v_{j'}\ov_{k'}+\ov_{j'}v_{k'}).
 \ea\ee
 By(\ref{ii33}), we get
 \bel{i33}\ba{ll}\ds
I_1^3 \overline {I_2^3}+\overline
 {I_1^3}I_2^3\\
\ns\ds=-2\sum_{j,k,j',k'=1}^m\[a^{jk}\ell_ja^{j'k'}(v_{j'}\ov_k+\ov_{j'}v_k)\]_{k'}+2\sum_{j,k,j',k'=1}^ma^{j'k'}(a^{jk}\ell_j)_{k'}(v_{j'}\ov_k+\ov_{j'}v_k)\\
\ns\ds\q+\sum_{j,k,j',k'=1}^m\[a^{jk}a^{j'k'}\ell_j(v_{j'}\ov_{k'}+\ov_{j'}v_{k'})\]_k-\sum_{j,k,j',k'=1}^m(a^{jk}a^{j'k'}\ell_j)_k(v_{j'}\ov_{k'}+\ov_{j'}v_{k'}).
 \ea\ee
Further,
 \bel{i34}\ba{ll}\ds
I_1^3 \overline {I_2^4}+\overline
 {I_1^3}I_2^4&\ds=\sum_{j,k=1}^m\[\Psi a^{jk}(v_j\ov+\ov_j
 v)\]_k-\sum_{j,k=1}^ma^{jk}\Psi(v_j\ov_k+\ov_jv_k)\\
 \ns&\ds\q-\sum_{j,k=1}^m\[a^{jk}\Psi_k|v|^2\]_j+\sum_{j,k=1}^m(a^{jk}\Psi_k)_j|v|^2.
 \ea\ee
Finally,
 \bel{i41}\ba{ll}\ds
I_1^4(\overline I_2^1+\overline I_2^2+\overline I_2^3+\overline
I_2^4)+\overline I_1^4
(I_2^1+I_2^2+I_2^3+I_2^4)\\
\ns\ds=(\a A|v|^2)_t-(\a A)_t|v|^2-2\sum_{j,k=1}^m(a^{jk}\ell_jA|v|^2)_k\\
\ns\ds\q+2\[\sum_{j,k=1}^m(a^{jk}\ell_jA)_k+A\Psi\]|v|^2.
 \ea\ee

\ms

{\it Step 3.} By (\ref{I})--(\ref{i41}), combining all
`$\frac{\pa}{\pa t}$-terms', all `$\frac{\pa}{\pa x_k}$-terms',
 and (\ref{1a8}), we arrive at the desired identity (\ref{2a2}).\endpf

 \ms

We have the following pointwise estimate for the complex parabolic
operator $Gz$.
\bc\label{qc1}  Let $b\in\dbR$ and $a^{jk}(t,x)\in
C^{1,2}(\dbR^{1+n};\;\dbR)$ satisfy condition (\ref{a1}). Let
$z,\;v\in C^2(\dbR^{1+m}; \;\dbC)$ and $\Psi,\;\ell\in
C^2(\dbR^{1+m};\dbR)$. Set $\th=e^\ell$ and $v=\th z$. Put
 \bel{Psi}
 \Psi=-2\sum_{j,k=1}^n(a^{jk}\ell_j)_k.
 \ee
Then
  \bel{lk5}\ba{ll}\ds
\th^2|G z|^2+M_t+\div V\\
\ns\ds\ge
\sum_{j,k,j',k'=1}^n\[2(a^{j'k}\ell_{j'})_{k'}a^{jk'}-a^{jk}_{k'}a^{j'k'}\ell_{j'}+\frac{1}{2}a^{jk}_t+
a^{jk}(a^{j'k'}\ell_{j'})_{k'}\](v_k\ov_j+\ov_k
v_j)\\
\ns\ds\qq+ib\sum_{j,k=1}^n(a^{jk}_t\ell_j+2a^{jk}\ell_{jt})(\ov_kv-v_k\ov)-ib\sum_{j,k=1}^n(a^{jk}\ell_j)_k(\ov
v_t-v\ov_t)+B|v|^2,
 \ea\ee
where
 \bel{c1a3}\left\{\ba{ll}\ds
M=\[(1+b^2)\ell_t-A\] |v|^2+\sum_{j,k=1}^na^{jk}v_j\ov_k+ib\sum_{j,k=1}^na^{jk}\ell_j(\ov_kv-v_k\ov),\\
 \ns\ds V^k=\sum_{j,j',k'=1}^n\Big\{-i b\[a^{jk}\ell_j(\ov_tv-\ov
 v_t)+a^{jk}\ell_t(v_j\ov-\ov_jv)\]-a^{jk}(v_j\ov_t+\ov_jv_t)\\
\ns\ds\qq\q +\(2a^{jk'}a^{j'k}-a^{jk}a^{j'k'}\)\ell_j(v_{j'}\ov_{k'}+\ov_{j'}v_{k'}) -\Psi a^{jk}(v_j\ov+\ov_jv)\\
 \ns\ds\qq\q
+a^{jk}(2A\ell_j+\Psi_j-2\ell_j\ell_t)|v|^2\Big\},\\
\ns\ds
B=(1+b^2)\ell_{tt}-A_t-2\[\sum_{j,k=1}^n(a^{jk}\ell_j\ell_{t})_k+\Psi\ell_t\]\\
\ns\ds\qq\q+\sum_{j,k=1}^n(a^{jk}\Psi_k)_j+2\[\sum_{j,k=1}^n(a^{jk}\ell_jA)_k+A\Psi\].
 \ea\right.\ee \ec

{\it Proof.} Recalling (\ref{cp2}) for the definition of $Gz$,
taking $ m=n, \;\a(x)\equiv 1,\;\b(x)\equiv b$ in Theorem \ref{2t1},
 by using H\"older inequality and simple computation, we immediately
 obtain
the desired result.\endpf

 \section{A modified point-wise estimate}\label{ss3}

Note that the term $\ds ib\sum_{j,k=1}^n(a^{jk}\ell_j)_k(\ov
v_t-v\ov_t)$ in the right-hand side of (\ref{lk5}) is not good. In
this section, we make some modification on this term and derive the
following modified point-wise inequality for the parabolic operator
with a complex principal.
\bt\label{qc2}  Let $b\in\dbR$ and $a^{jk}(t,x)\in
C^{1,2}(\dbR^{1+n};\;\dbR)$ satisfy (\ref{a1})--(\ref{a2}). Let
$z,\;v\in C^2(\dbR^{1+m}; \;\dbC)$ and $\Psi,\;\ell\in
C^2(\dbR^{1+m};\dbR)$ satisfy (\ref{Psi}). Put
 \bel{gpoint}
 \th=e^\ell,\q v=\th z.
 \ee
Then
  \bel{qsds}\ba{ll}\ds
2\th^2|Gz|^2+M_t+\div \tilde
V\\
\ns\ds\ge
\sum_{j,k=1}^nc^{jk}(v_k\ov_j+\ov_kv_j)+\tilde B|v|^2\\
 \ns\ds\q+ib\sum_{j,k=1}^n\[a^{jk}_t\ell_j+2a^{jk}\ell_{jt}+{1\over
{1+b^2}}\sum_{j',k'=1}^n(a^{j'k'}\ell_{j'})_{k'j}a^{jk}\](\overline
v_kv-v_k\overline v),
 \ea\ee
where $M,\;V^k,\; B$ is given by (\ref{c1a3}) and
 \bel{qagdf3}\left\{\ba{ll}\ds
  \tilde
V^k=V^k+\sum_{j,j',k'=1}^n\Big\{{ib\over
{1+b^2}}\sum_{j,k,j',k'=1}^n\[(a^{j'k'}\ell_{j'})_{k'}a^{jk}(\overline
v_jv-v_j\overline v)\]\\
\ns\ds\qq\q-{b^2\over
{1+b^2}}\th^2(a^{j'k'}\ell_{j'})_{k'}a^{jk}(\overline
z_jz+z_j\overline z)+{b^2\over
{1+b^2}}\[\th^2(a^{j'k'}\ell_{j'})_{k'}\]_ja^{jk}|z|^2\\
\ns\ds\qq\q-{2b^2\over {1+b^2}}
(a^{j'k'}\ell_{j'})_{k'}a^{jk}\ell_j|v|^2\Big\},\\
 \ns\ds
c^{jk}=\sum_{j',k'=1}^n\[2(a^{j'k}\ell_{j'})_{k'}a^{jk'}-a^{jk}_{k'}a^{j'k'}\ell_{j'}+\frac{1}{2}a^{jk}_t+{1\over
{2(1+b^2)}}
a^{jk}(a^{j'k'}\ell_{j'})_{k'}\],\\
 \ns\ds \tilde B=B-{2b^2\over
{1+b^2}}\sum_{j,k,j',k'=1}^n(a^{j'k'}\ell_{j'})_{k'}a^{jk}\ell_k\ell_j-{{b^2}\over{1+b^2}}\|\sum_{j,k=1}^n(a^{jk}\ell_j)_k\|^2\\
 \ns\ds\qq\q+{b^2\over
{1+b^2}}\sum_{j,k,j',k'=1}^n\Big\{2a^{jk}\ell_k(a^{j'k'}\ell_{j'})_{k'j}+\(a^{jk}(a^{j'k'}\ell_{j'})_{k'j}\)_k\Big\}.
 \ea\right.\ee
\et

{\it Proof.} We divided the proof into several steps.

\ms

{\it Step 1. } Note that $v=\th z$, it is easy to check that
 \bel{ww4}
\ov v_t-v\ov_t=\th^2 (\overline z z_t-z{\overline z}_t),\q \ov
v_k-v\ov_k=\th^2 (\overline z z_k-z{\overline z}_k).
 \ee
Recalling (\ref{cp2}) for the definition of $Gz$ and by (\ref{ww4}),
we have
 \bel{ib2}\ba{ll}&\ds
 -ib\sum_{j,k=1}^n(a^{jk}\ell_j)_k(\ov
v_t-v\ov_t)=-ib\sum_{j,k=1}^n(a^{jk}\ell_j)_k\th^2(\overline{z}
z_t-z\overline{z}_t)\\
\ns&\ds=-{ib\th^2\over
{1+b^2}}\sum_{j,k=1}^n(a^{jk}\ell_j)_k\[(1-ib)G z\overline
z-\overline {(1-ib) G z} z\]\\
\ns&\ds\q -{ib\th^2\over
{1+b^2}}\sum_{j',k'=1}^n(a^{j'k'}\ell_{j'})_{k'}\sum_{j,k=1}^n\[(a^{jk}\overline
z_j)_kz-(a^{jk}z_j)_k\overline z)\]\\
\ns&\ds\q+{b^2\th^2\over
{1+b^2}}\sum_{j',k'=1}^n(a^{j'k'}\ell_{j'})_{k'}\sum_{j,k=1}^n\[(a^{jk}\overline
z_j)_kz+(a^{jk}z_j)_k\overline z)\]\\
\ns&\ds\equiv\sum_{k=1}^3J_k.
 \ea\ee

\ms

{\it Step 2. } Let us estimate $J_k\;(k=1,2,3)$ respectively.

First, note that $v=\th z$, we have
 \bel{kl1}\ba{ll}\ds
J_1&\ds= -{ib\th^2\over
{1+b^2}}\sum_{j,k=1}^n(a^{jk}\ell_j)_k\[(1-ib)G z\overline
z-\overline {(1-ib) G z} z\]\\
 \ns&\ds\ge -\|{{(1-ib)\th
G z}\over{\sqrt{1+b^2}}}\|^2-\|{{ib\th
z}\over{\sqrt{1+b^2}}}\sum_{j,k=1}^n(a^{jk}\ell_j)_k\|^2\\
\ns&\ds=-\th^2| G
z|^2-{{b^2}\over{1+b^2}}\|\sum_{j,k=1}^n(a^{jk}\ell_j)_k\|^2|v|^2.
 \ea\ee
Next, by using (\ref{ww4}) again, we have
 \bel{ww5}\ba{ll}\ds
J_2&\ds=-{ib\th^2\over
{1+b^2}}\sum_{j',k'=1}^n(a^{j'k'}\ell_{j'})_{k'}\sum_{j,k=1}^n\[(a^{jk}\overline
z_j)_kz-(a^{jk}z_j)_k\overline z)\]\\
\ns&\ds =-{ib\over
{1+b^2}}\sum_{j,k,j',k'=1}^n\[\th^2(a^{j'k'}\ell_{j'})_{k'}a^{jk}(\overline
z_jz-z_j\overline z)\]_k\\
\ns&\ds\q+{ib\over
{1+b^2}}\sum_{j,k,j',k'=1}^n\th^2\[2\ell_k(a^{j'k'}\ell_{j'})_{k'}+(a^{j'k'}\ell_{j'})_{k'k}\]a^{jk}(\overline
z_jz-z_j\overline z)\\
\ns&\ds =-{ib\over
{1+b^2}}\sum_{j,k,j',k'=1}^n\[(a^{j'k'}\ell_{j'})_{k'}a^{jk}(\overline
v_jv-v_j\overline v)\]_k\\
\ns&\ds\q+{ib\over
{1+b^2}}\sum_{j,k,j',k'=1}^n\[2\ell_k(a^{j'k'}\ell_{j'})_{k'}+(a^{j'k'}\ell_{j'})_{k'k}\]a^{jk}(\overline
v_jv-v_j\overline v).
 \ea\ee
Next, noting that $a^{jk}$ satisfy (\ref{a1}), and recalling $v=\th
z$, it follows
 \bel{kl2}\ba{ll}\ds
J_3&\ds={b^2\th^2\over
{1+b^2}}\sum_{j',k'=1}^n(a^{j'k'}\ell_{j'})_{k'}\sum_{j,k=1}^n\[(a^{jk}\overline
z_j)_kz+(a^{jk}z_j)_k\overline z)\]\\
\ns&\ds ={b^2\over
{1+b^2}}\sum_{j,k,j',k'=1}^n\Big\{\[\th^2(a^{j'k'}\ell_{j'})_{k'}a^{jk}(\overline
z_jz+z_j\overline z)\]-\[\th^2(a^{j'k'}\ell_{j'})_{k'}\]_ja^{jk}|z|^2\Big\}_k\\
\ns&\ds\q-{b^2\th^2\over
{1+b^2}}\sum_{j,k,j',k'=1}^n(a^{j'k'}\ell_{j'})_{k'}a^{jk}(\overline
z_jz_k+z_j\overline z_k)\\
\ns&\ds\q+{b^2\over
{1+b^2}}\sum_{j,k,j',k'=1}^n\Big\{\[\th^2(a^{j'k'}\ell_{j'})_{k'}\]_ja^{jk}\Big\}_k|z|^2\\
\ns&\ds={b^2\over
{1+b^2}}\sum_{j,k,j',k'=1}^n\Big\{\[\th^2(a^{j'k'}\ell_{j'})_{k'}a^{jk}(\overline
z_jz+z_j\overline z)\]-\[\th^2(a^{j'k'}\ell_{j'})_{k'}\]_ja^{jk}|z|^2\Big\}_k\\
\ns&\ds\q+{2b^2\over
{1+b^2}}\sum_{j,k,j',k'=1}^n\[(a^{j'k'}\ell_{j'})_{k'}a^{jk}\ell_k|v|^2\]_j\\
\ns&\ds\q-{b^2\over
{1+b^2}}\sum_{j,k,j',k'=1}^n(a^{j'k'}\ell_{j'})_{k'}a^{jk}(\overline
v_jv_k+v_j\overline
v_k)+\Big\{2\sum_{j,k,j',k'=1}^na^{jk}\ell_j\ell_k(a^{j'k'}\ell_{j'})_{k'}\\
\ns&\ds\qq\qq\q+2\sum_{j,k,j',k'=1}^na^{jk}\ell_k(a^{j'k'}\ell_{j'})_{k'j}+\sum_{j,k,j',k'=1}^n\(a^{jk}(a^{j'k'}\ell_{j'})_{k'j}\)_k\Big\}|v|^2
 \ea\ee
where we have used the following fact:
 \bel{kll3}\ba{ll}\ds
\th^2\sum_{j,k,j',k'=1}^n(a^{j'k'}\ell_{j'})_{k'}a^{jk}(\overline
z_jz_k+z_j\overline z_k)\\
\ns\ds=\sum_{j,k,j',k'=1}^n(a^{j'k'}\ell_{j'})_{k'}a^{jk}(\overline
v_jv_k+v_j\overline
v_k)+2\sum_{j,k,j',k'=1}^n(a^{j'k'}\ell_{j'})_{k'}a^{jk}\ell_j\ell_k|v|^2\\
\ns\ds\q-2\sum_{j,k,j',k'=1}^n\[(a^{j'k'}\ell_{j'})_{k'}a^{jk}\ell_k|v|^2\]_j+2\sum_{j,k,j',k'=1}^n\[(a^{j'k'}\ell_{j'})_{k'}a^{jk}\ell_k\]_j|v|^2.
 \ea\ee
and
 \bel{klk3}\ba{ll}\ds
\sum_{j,k,j',k'=1}^n\Big\{\[\th^2(a^{j'k'}\ell_{j'})_{k'}\]_ja^{jk}\Big\}_k|z|^2\\
\ns\ds=\sum_{j,k,j',k'=1}^n\Big\{\th^2\[2a^{jk}\ell_j(a^{j'k'}\ell_{j'})_{k'}+a^{jk}(a^{j'k'}\ell_{j'})_{k'j}\]\Big\}_k|z|^2\\
\ns\ds=\Big\{4\sum_{j,k,j',k'=1}^na^{jk}\ell_j\ell_k(a^{j'k'}\ell_{j'})_{k'}+4\sum_{j,k,j',k'=1}^na^{jk}\ell_k(a^{j'k'}\ell_{j'})_{k'j}\\
\ns\ds\qq\q+2\|\sum_{j,k=1}^n(a^{jk}\ell_j)_k\|^2+\sum_{j,k,j',k'=1}^n\(a^{jk}(a^{j'k'}\ell_{j'})_{k'j}\)_k\Big\}|v|^2.
 \ea\ee

\ms

{\it Step 3. } Noting that $a^{jk}$ satisfying (\ref{a1})--
(\ref{a2}), we have the following fact.
 \bel{kkll3}\ba{ll}\ds
{1\over2(1+b^2)}\sum_{j,k=1}^na^{jk}(\overline v_jv_k+v_j\overline
v_k)-2i{b\over{1+b^2}}\sum_{j,k=1}^na^{jk}\ell_k(\overline
v_jv-v_j\overline
v)\\
\ns\ds={1\over{1+b^2}}\sum_{j,k=1}^na^{jk}\[v_j\overline
v_k-2ib\ell_k(\overline v_jv-v_j\overline
v)+4b^2\ell_j\ell_k|v|^2\]-{4b^2\over{1+b^2}}\sum_{j,k=1}^na^{jk}\ell_j\ell_k|v|^2\\
\ns\ds={1\over{1+b^2}}\sum_{j,k=1}^na^{jk}(v_j-2ib\ell_jv)\overline{(v_k-2ib\ell_kv)}-{4b^2\over{1+b^2}}\sum_{j,k=1}^na^{jk}\ell_j\ell_k|v|^2\\
\ns\ds\ge-{4b^2\over{1+b^2}}\sum_{j,k=1}^na^{jk}\ell_j\ell_k|v|^2.
 \ea\ee

Finally, combining (\ref{lk5}), (\ref{ww4})--(\ref{kl2}) and
(\ref{kkll3}), we arrive at the desired result (\ref{qsds}).\endpf

\section{Global Carleman estimate for parabolic operators with complex principal part}\label{ss4}

We begin with the following known result.

\bl\label{h31} {\rm (\cite{FI}, \cite{WW})} There is a real function
$\psi\in C^2(\oO)$ such that $\psi>0$ in $\O$ and $\psi=0$ on
$\pa\O$ and $|\n\psi(x)|>0$ for all $x\in\oO\setminus\o_0$.
 \el

 For any (large) parameters $\l>1$ and $\mu>1$, put

 \bel{p33}
 \ell=\l\rho,\q\varphi(t,x)={e^{\mu\psi(x)}\over{t(T-t)}},\q\rho(t,x)={{e^{\mu\psi(x)}-e^{2\mu|\psi|_{C(\cl{\O
 };\;\dbR)}}}\over{t(T-t)}}.
 \ee

For $j,k=1,\cdots,n$, it is easy to check that
 \bel{33ss2}
 \ell_t=\l\rho_t,\q \ell_j=\l\mu\varphi\psi_j,\q
 \ell_{jk}=\l\mu^2\varphi\psi_j\psi_k+\l\mu\varphi\psi_{jk}\ee
 and
  \bel{3a3ss2}
  |\rho_t|\le C e^{2\mu |\psi|_{C(\overline\O)}}\varphi^2,\q |\varphi_t|\le
 C\varphi^2.
 \ee
In the sequel, for $k\in\dbN$, we denote by $O(\mu^k)$ a function of
order $\mu^k$ for large $\mu$ (which is independent of $\l$); by
$O_\mu(\l^k)$ a function of order $\l^k$ for fixed $\mu$ and for
large $\l$.

We have the following Carleman estimate for the differential
operator $G$ defined in (\ref{cp2}):

 \bt\label{gt3} Let $b\in\dbR$ and $a^{jk}$ satisfy (\ref{a1})--(\ref{a2}).
Then there is a $\mu_0>0$ such that for all $\mu\ge\mu_0$, one can
find two constants $C=C(\mu)>0$ and $\l_1=\l_1(\mu)$, such that for
all  $z\in C([0,T]; L^2(\O))\bigcap C((0,T]; H_0^1(\O))$ and  for
all $\l\ge\l_1$, it holds
 \bel{gcpp}\ba{ll}\ds
 \l^3\mu^4\int_Q\varphi^3\th^2|z|^2dtdx+\l\mu^2\int_Q\varphi\th^2|\n
 z|^2dtdx\\
 \ns\ds\q\le
 C(1+b^2
)\[\int_Q\th^2|G
z|^2dtdx+\l^3\mu^4\int_{(0,T)\t\o}\!\varphi^3\th^2|z|^2dtdx\].
 \ea\ee \et

\ms

{\it Proof.} The proof is long, we divided it into several steps.

\ms

{\it Step 1.} Recalling (\ref{qagdf3}) for the definition of
$c^{jk}$, by (\ref{33ss2})--(\ref{3a3ss2}), we have
 \bel{cjk}\ba{ll}\ds\sum_{j,k=1}^nc^{jk}(v_k\ov_j+\ov_kv_j)\\
 \ns\ds=\sum_{j,k,j',k'=1}^n\[2(a^{j'k}\ell_{j'})_{k'}a^{jk'}-a^{jk}_{k'}a^{j'k'}\ell_{j'}\\
 \ns\ds\qq\qq\qq\qq\qq\qq\q+\frac{1}{2}a^{jk}_t+{1\over{2(1+b^2)}}
a^{jk}(a^{j'k'}\ell_{j'})_{k'}\](v_k\ov_j+\ov_kv_j)\\
 \ns\ds=4\l\mu^2\varphi\|\sum_{j,k=1}^na^{jk}\psi_j\ov_k\|^2\\
 \ns\ds\qq+{1\over{2(1+b^2)}}\l\mu^2\varphi\sum_{j,k,j',k'=1}^na^{jk}a^{j'k'}\psi_{j'}\psi_{k'}(v_k\ov_j+\ov_kv_j)+\l\varphi
 O(\mu)|\n v|^2.
 \ea\ee
On the other hand, by (\ref{33ss2})--(\ref{3a3ss2}), recalling
(\ref{Psi}) and (\ref{1f3}) for the definitions of $\Psi$ and $A$,
it is easy to check that
 \bel{ghass1}
\Psi=-2\l\mu^2\varphi\sum_{j,k=1}^na^{jk}\psi_j\psi_k+\l\varphi
O(\mu),\q
A=\l^2\mu^2\varphi^2\sum_{j,k=1}^na^{jk}\psi_j\psi_k+\varphi^2
O_\mu(\l).
 \ee
Next, recalling (\ref{1f3}) and (\ref{c1a3}) for the definition of
$A$ and $B$, respectively.  By (\ref{33ss2})--(\ref{3a3ss2}) and
combining (\ref{ghass1}), we have
 \bel{gbhb}\ba{ll}\ds
 B&\ds=2\sum_{j,k=1}^na^{jk}\ell_jA_k-2 A\sum_{j,k=1}^n(a^{jk}\ell_j)_k\\
\ns&\ds\qq+(1+b^2)\ell_{tt}-A_t-2\sum_{j,k=1}^na^{jk}\ell_j\ell_{tk}+2\ell_t\sum_{j,k=1}^n(a^{jk}\ell_j)_k+\sum_{j,k=1}^n(a^{jk}\Psi_k)_j\\
 \ns&\ds
 =2\l^3\mu^4\varphi^3\|\sum_{j,k=1}^na^{jk}\psi_j\psi_k\|^2+\l^3\varphi^3O(\mu^3)+\varphi^3O_\mu(\l^2).
 \ea\ee
Hence, by recalling (\ref{qagdf3}) for the definition of $\tilde B$,
we have
 \bel{kl4}\ba{ll}\ds
\tilde
B&\ds={2\over{1+b^2}}\l^3\mu^4\varphi^3\|\sum_{j,k=1}^na^{jk}\psi_j\psi_k\|^2+\l^3\varphi^3O(\mu^3)+\varphi^3O_\mu(\l^2).
 \ea\ee

{\it Step 2.} By (\ref{33ss2})--(\ref{3a3ss2}), we have
 \bel{op4}\ba{ll}\ds
\|ib\sum_{j,k=1}^n(a^{jk}_t\ell_j+2a^{jk}\ell_{jt})(\ov_kv-v_k\ov)\|\\
\ns\ds\le C\l\mu\|\varphi^2\sum_{j,k=1}^na^{jk}\psi_j(\overline
v_kv-v_k\overline v)\|\le
C\l\mu\varphi\|\sum_{j,k=1}^na^{jk}\psi_j\ov_k\|^2+C\l\mu\varphi^3|v|^2.
 \ea
 \ee
It is easy to see that (\ref{op4}) can be absorbed by (\ref{cjk})
and (\ref{kl4}).

Similarly, by using (\ref{33ss2})--(\ref{3a3ss2}) again, we have
 \bel{lk7}\ba{ll}\ds
 \|{ib\over
{1+b^2}}\sum_{j,k,j',k'=1}^n(a^{j'k'}\ell_{j'})_{k'j}a^{jk}(\overline
v_kv-v_k\overline v)
\|\\
\ns\ds\le {1\over
{1+b^2}}\|b\[\l\mu^3\varphi^2\sum_{j',k'=1}^na^{j'k'}\psi_{j'}\psi_{k'}+\l\varphi^2O(\mu^2)\]a^{jk}\psi_j(\overline
v_kv-v_k\overline v)\|\\
\ns\ds\le {1\over
{1+b^2}}\|b\l\mu^3\varphi^2\sum_{j',k'=1}^na^{j'k'}\psi_{j'}\psi_{k'}a^{jk}\psi_j(\overline
v_kv-v_k\overline
v)\|\\
\ns\ds\q+{C\over
{1+b^2}}\|b\l\mu^2\varphi^2\sum_{j,k=1}^na^{jk}\psi_j(\overline
v_kv-v_k\overline v)\|\\
\ns\ds\le  {1 \over
{1+b^2}}\l\mu^2\varphi\|\sum_{j,k=1}^na^{jk}\psi_j\overline
v_k\|^2+{ b^2\over
{1+b^2}}\l\mu^4\varphi^3\|\sum_{j,k=1}^na^{jk}\psi_j\psi_k\|^2|v|^2\\
\ns\ds\qq+{C\over
{1+b^2}}\l\mu\varphi\|\sum_{j,k=1}^na^{jk}\psi_j\ov_k\|^2+{C
b^2\over {1+b^2}}\l\mu^3\varphi^3|v|^2.
 \ea\ee
Therefore (\ref{lk7}) also can be absorbed by (\ref{cjk}) and
(\ref{kl4}).

Combining (\ref{cjk}), (\ref{kl4})--(\ref{lk7}), by (\ref{a2}), we
end up with
 \bel{lk9}\ba{ll}\ds
\hbox{ The right-hand side of (\ref{qsds})}\\
\ns\ds \ge {1\over
{1+b^2}}\l\mu^2\varphi\sum_{j',k'=1}^na^{j'k'}\psi_{j'}\psi_{k'}\sum_{j,k=1}^na^{jk}v_k\ov_j+\l\varphi
 O(\mu)|\n v|^2\\
 \ns\ds\q+\[{2\over{1+b^2}}\l^3\mu^4\varphi^3\|\sum_{j,k=1}^na^{jk}\psi_j\psi_k\|^2+\l^3\varphi^3O(\mu^3)+\varphi^3O_\mu(\l^2)\]|v|^2\\
 \ns\ds \ge
{1\over{1+b^2}}\[s_0^2\l\mu^2\varphi|\n\psi|^2+\l\varphi
 O(\mu)\]|\n v|^2\\
 \ns\ds\q+{2\over{1+b^2}}\[s_0^2\l^3\mu^4\varphi^3|\n\psi|^4+\l^3\varphi^3O(\mu^3)+\varphi^3O_\mu(\l^2)\]|v|^2.
 \ea\ee

 \ms

{\it Step 3.} Integrating (\ref{qsds}) on $Q$,  by (\ref{lk9}),
noting that $\th (0)=\th (T)\equiv 0$, we have
 \bel{lk10}\ba{ll}\ds
\int_Q\varphi\[\l\mu^2|\n\psi|^2+\l
 O(\mu)\]|\n v|^2dxdt\\
 \ns\ds\q+\int_Q\varphi^3\[\l^3\mu^4|\n\psi|^4+\l^3O(\mu^3)+O_\mu(\l^2)\]|v|^2dxdt\\
 \ns\ds\le C\[|\th Gz|^2_{L^2(Q)}+\int_Q \div \tilde V\cd\nu dxdt\].
 \ea\ee
Next, recalling (\ref{1a3}) and (\ref{qagdf3}) for the definition of
$V$ and $\tilde V$.  Noting that $z|_\Si=0$ and $ v_i={{\pa v}\over
{\pa\nu}}\nu_i$ (which follows from (\ref{gap1}) and $v|_\Si=0$,
respectively), by (\ref{33ss2}) and Lemma \ref{h31}, we have
 \bel{bv}\ba{ll}\ds
\int_Q \div \tilde V\cd\nu dxdt=\int_Q \div V\cd\nu dxdt\\
 \ns\ds=\l\mu\int_\Si\varphi {{\pa\psi}\over {\pa\nu}}\|{{\pa
v}\over {\pa\nu}}\|^2\(\sum_{j,k=1}^na^{jk}\nu_j\nu_k\)^2dtdx \le 0.
 \ea\ee

On the other hand,
 \bel{lk11}\ba{ll}\ds
\hb {Left-side of
(\ref{lk10})}&\ds=\int_0^T\(\int_{\O\setminus{\o_0}}+\int_{\o_0}\)
\[\varphi \(\l\mu^2|\n \psi|^2+\l O(\mu)\)|\n
  v|^2\\
  \ns&\ds\qq+\varphi^3\(\l^3\mu^4 |\n\psi|^4+\l^3
O(\mu^3)+\varphi^3 O_\mu(\l^2)\)|v|^2 \]dtdx\\
\ns&\ds\ge \int_0^T\int_{\O\setminus{\o_0}}\[\varphi \(\l\mu^2|\n
\psi|^2+\l O(\mu)\)|\n
  v|^2\\
  \ns&\ds\qq+\varphi^3\(\l^3\mu^4 |\n\psi|^4+\l^3
O(\mu^3)+\varphi^3 O_\mu(\l^2)\)|v|^2 \]dtdx\\
\ns&\ds\qq-C\l\mu^2\int_{Q_0}\varphi(|\n v|^2+\l^2\mu^2\varphi^2
|v|^2)dtdx. \ea\ee

By the choose of $\psi$, we know that $\ds\min_{x\in {\O\setminus
{\o_0}}}|\n \psi|>0$. Hence, there is a $\mu_0>0$ such that for all
$\mu\ge\mu_0$, one can find a constant $\l_1=\l_1(\mu)$ so that for
any $\l\ge\l_1$, it holds
 \bel{kpa4s6}
\ba{ll}\ds \int_0^T\int_{\O\setminus{\o_0}} \varphi \[\l\mu^2|\n
\psi|^2+\l O(\mu)\]|\n
  v|^2 dtdx\\
  \ns\ds\qq+\int_Q \varphi^3\[\l^3\mu^4 |\n\psi|^4+\l^3
O(\mu^3)+\varphi^3 O_\mu(\l^2)\]|v|^2 dtdx\\
\ns\ds\ge c_0\l\mu^2\int_0^T\int_{\O\setminus{\o_0}}\varphi(|\n
v|^2+\l^2\mu^2\varphi^2 |v|^2)dtdx, \ea\ee
 where $\ds c_0=\min\(\min_{x\in {\O\setminus
{\o_0}}}|\n \psi|^2,\min_{x\in {\O\setminus {\o_0}}}|\n \psi|^4\)$
is a positive constant.

Next, note that $v=\th z$, by (\ref{33ss2}), we have
 \bel{bb1}
z_j=\th^{-1}(v_j-\ell_jv)=\th^{-1}(v_j-\l\mu\varphi\psi_j v),\q
v_j=\th(z_j+\ell_jz)=\th(z_j+\l\mu\varphi\psi_j z).
 \ee
By (\ref{bb1}), we get
 \bel{bb2}
{1\over C}\th^2(|\n z|^2+\l^2\mu^2\varphi^2 |z|^2)\le|\n
v|^2+\l^2\mu^2\varphi^2 |v|^2\le C\th^2(|\n z|^2+\l^2\mu^2\varphi^2
|z|^2).
 \ee
Therefore, it follows from (\ref{lk11})--(\ref{kpa4s6})  and
(\ref{bb2}) that
 \bel{kpa4s7}\ba{ll}\ds
 \l\mu^2\int_Q\varphi\th^2(|\n z|^2+\l^2\mu^2\varphi^2 |z|^2)dtdx\\
 \ns\ds=\l\mu^2\int_0^T\(\int_{\O\setminus{\o_0}}+\int_{\o_0}\) \varphi\th^2(|\n z|^2+\l^2\mu^2\varphi^2 |z|^2)dtdx\\
\ns\ds\le C\Big\{ \int_Q \varphi \[\l\mu^2|\n \psi|^2+\l O(\mu)\]|\n
  v|^2 dtdx\\
  \ns\ds\qq\q+\int_Q \varphi^3\[\l^3\mu^4 |\n\psi|^4+\l^3
O(\mu^3)+\varphi^3 O_\mu(\l^2)\]|v|^2 dtdx\\
\ns\ds\qq\q+\l\mu^2\int_{Q_0} \varphi\th^2(|\n
z|^2+\l^2\mu^2\varphi^2 |z|^2)dtdx\Big\}. \ea\ee

Now, combining (\ref{lk10})--(\ref{bv}) and (\ref{kpa4s7}), we end
up with
 \bel{kg4s8}\ba{ll}\ds
\l\mu^2\int_Q\varphi\th^2(|\n z|^2+\l^2\mu^2\varphi^2 |z|^2)dtdx\\
\ns\ds\le C\[\int_Q\th^2|Gz|^2dtdx+\l\mu^2\int_{Q_0}
\varphi\th^2(|\n z|^2+\l^2\mu^2\varphi^2 |z|^2)dtdx\]. \ea\ee

\ms

{\it Step 4.} Finally, we choose a cut-off function $\zeta\in
C_0^\i(\o;\;[0,1])$ so that $\zeta\equiv 1$ on $\o_0$. Then
 \bel{bb3}
 \[\zeta^2\varphi\th^2|(1+ib)z|^2\]_t=\zeta^2|(1+ib)z|^2(\varphi\th^2)_t+\zeta^2\varphi\th^2(1+ib)(1-ib)(\overline zz_t+z\overline
 z_t).
 \ee
By (\ref{cp2}), (\ref{bb3}) and noting $\ds\th(0,x)=\th(T,x)\equiv
0$, we find
 \bel{gg4s9}\ba{ll}
0&\ds=\int_{Q_0}\zeta^2\[|(1+ib)z|^2(\varphi\th^2)_t+\varphi\th^2(1+ib)(1-ib)(\overline
zz_t+z\overline z_t)\]dtdx\\
\ns&\ds=\int_{Q_0}\zeta^2\th^2\Big\{|(1+ib)z|^2(\varphi_t+2\l\varphi\rho_t)+\varphi(1-ib)\overline
z\[-\sum_{j,k}(a^{jk}z_j)_k+G z\]\\
\ns&\ds\qq\q+\varphi(1+ib) z\[-\sum_{j,k}(a^{jk}\overline
z_j)_k+\overline {G z}\]\Big\}\\
\ns&\ds=\int_{Q_0}\th^2\Big\{\zeta^2|(1+ib)z|^2(\varphi_t+2\l\varphi\rho_t)\\
\ns&\ds\qq+\zeta^2\varphi \sum_{j,k}a^{jk}\[(1-ib)z_j\overline z_k+(1+ib)\overline z_j z_k\]\\
\ns&\ds\qq+\mu\zeta^2\varphi\sum_{j,k}a^{jk}\[(1-ib)\overline
zz_j\psi_k+(1+ib) z\overline
z_j\psi_k\]\\
\ns&\ds\qq+2\l\mu\zeta^2\varphi^2\sum_{j,k}a^{jk}\[(1-ib)\overline
zz_j\psi_k+(1+ib)
z\overline z_j\psi_k\]\\
\ns&\ds\qq +2\zeta\varphi \sum_{j,k}a^{jk}\[(1-ib)\overline
zz_j\zeta_k+(1+ib)z\overline z_j\zeta_k\]\\
\ns&\ds\qq+\zeta^2\varphi [(1-ib)\overline z G z+(1+ib)z\overline {G
z}]\Big\}dtdx. \ea\ee

Hence, by (\ref{gg4s9}) and (\ref{a2}), we conclude that, for some
$\d>0$,
 \bel{gg4s10}\ba{ll}\ds
 2\int_{Q_0}\zeta^2\varphi\th^2|\n z|^2dxdt\\
\ns\ds=\|\int_{Q_0}\th^2\Big\{\zeta^2|(1+ib)z|^2(\varphi_t+2\l\varphi\rho_t)\\
\ns\ds\qq+\mu\zeta^2\varphi\sum_{j,k}a^{jk}\[(1-ib)\overline
zz_j\psi_k+(1+ib) z\overline
z_j\psi_k\]\\
\ns\ds\qq+2\l\mu\zeta^2\varphi^2\sum_{j,k}a^{jk}\[(1-ib)\overline
zz_j\psi_k+(1+ib)
z\overline z_j\psi_k\]\\
\ns\ds\qq +2\zeta\varphi \sum_{j,k}a^{jk}\[(1-ib)\overline
zz_j\zeta_k+(1+ib)z\overline z_j\zeta_k\]\\
\ns\ds\qq+\zeta^2\varphi [(1-ib)\overline z G z+(1+ib)z\overline
{G z}]\Big\}dtdx\|\\
\ns\ds\le\d\int_{Q_0}\zeta^2\varphi\th^2|\n z|^2dtdx\\
\ns\ds\qq+{C\over\d}\[{(1+b^2)\over{\l^2\mu^2}}\int_{Q_0}\th^2|G
z|^2dtdx+\l^2\mu^2\int_{Q_0}\varphi^3\th^2|z|^2dtdx\].
 \ea\ee
Now, we choose $\d=1$. By (\ref{gg4s10}), we conclude that
 \bel{gg4s11}
\int_{Q_0}\varphi\th^2|\n z|^2dtdx\le
C(1+b^2)\[{1\over{\l^2\mu^2}}\int_{Q_0}\th^2|Gz|^2dtdx+\l^2\mu^2\int_{Q_0}\varphi^3\th^2|z|^2dtdx\].
 \ee
Finally, combining (\ref{kg4s8}) and (\ref{gg4s11}), we arrive
\bel{gpp}\ba{ll}\ds
 \l^3\mu^4\int_Q\varphi^3\th^2|z|^2dtdx+\l\mu^2\int_Q\varphi\th^2|\n
 z|^2dtdx\\
 \ns\ds\le
 C(1+b^2)\[\int_Q\th^2|G z|^2dtdx+\l^3\mu^4\int_{(0,T)\t\o}\!\varphi^3\th^2|z|^2dtdx\],
 \ea\ee
which gives the proof of Theorem \ref{gt3}.\endpf

 \section {Proof of the main results}\label{ss5}

In this section, we will give the proofs of Theorem \ref{gt5} and
\ref{gt11}. Thanks to the classical dual argument and the fixed
point technique, proceeding as \cite{DFGZ}, Theorem \ref{gt11} is a
consequence of the observability result (\ref{3g1}). Therefore, we
only give here a brief proof of Theorem \ref{gt5}.

{\it Proof of Theorem \ref{gt5}.} We apply Theorem \ref{gt3} to
system (\ref{gap1}). Recalling that $q(\cd)\in L^\i(0,T;L^n(\O))$
and using (\ref{33ss2}), we obtain
 \bel{ed1}\ba{ll}\ds
\l^3\mu^4\int_Q\varphi^3\th^2|z|^2dtdx+\l\mu^2\int_Q\varphi\th^2|\n
 z|^2dtdx\\
 \ns\ds\le
 C(1+b^2)\[\int_Q\th^2|qz
|^2dtdx+\l^3\mu^4\int_{(0,T)\t\o}\!\varphi^3\th^2|z|^2dtdx\]\\
\ns\ds\le   C(1+b^2)(1+|r|^2)\[\int_Q\th^2|z
|^2dtdx+\l^3\mu^4\int_{(0,T)\t\o}\!\varphi^3\th^2|z|^2dtdx\].
 \ea\ee
Choosing $\ds\l\ge  C(1+b^2)(1+|r|^2)$ and $\mu$ large enough, from
(\ref{ed1}), one deduces that
 \bel{ed2}
\int_Q\varphi^3\th^2|z|^2dtdx\le
C\int_{(0,T)\t\o}\!\varphi^3\th^2|z|^2dtdx.
 \ee
Finally, by (\ref{ed2}) and  applying the usual energy estimate to
system (\ref{gap1}), we conclude that inequality (\ref{3g1}) holds,
with the observability constant $\cC$ given by (\ref{ccg}), which
completes the proof of Theorem \ref{gt5}.\endpf

\end{document}